\documentclass[11pt]{article}

\title{\bf Selected results on selection principles}

\author{\Large Ljubi\v sa D.R. Ko\v cinac }

\date{}

\usepackage{amsfonts}
\usepackage{amssymb}

\newtheorem{theorem}{{\bf Theorem}}

\newtheorem{problem}[theorem]{{\bf Problem}}

\newtheorem{conjecture}[theorem]{{\bf Conjecture}}

\newcommand{\sone}{{\sf S}_1}
\newcommand{\gone}{{\sf G}_1}
\newcommand{\sfin}{{\sf S}_{fin}}
\newcommand{\gfin}{{\sf G}_{fin}}

\newcommand{\ssone}{{\sf S}_1^*}

\newcommand{\ssfin}{{\sf S}_{fin}^*}

\newcommand{\sssone}{{\sf SS}_1^*}

\newcommand{\sssfin}{{\sf SS}_{fin}^*}

\newcommand{\sufin}{{\sf U}_{fin}^*}

\newcommand{\ckx}{{\sf C}_k(X)}
\newcommand{\cpx}{{\sf C}_p(X)}
\newcommand{\cpy}{{\sf C}_p(Y)}
\newcommand{\zero}{\underline{\sf 0}}

\newcommand{\fplusx}{(2^X,{\sf F}^+)}
\newcommand{\zplusx}{(2^X,{\sf Z}^+)}

\newcommand{\naturals}{{\mathbb N}}
\newcommand{\reals}{{\mathbb R}}

\newcommand{\integers}{{\mathbb Z}}

\begin{document}
\maketitle

\begin{abstract}
We review some selected recent results concerning selection
principles in topology and their relations with several
topological constructions.
\end{abstract}

\begin{flushleft}
\vspace{.2cm} %
\end{flushleft}

\section{Introduction}

The beginning of investigation of covering properties of
topological spaces defined in terms of diagonalization and
nowadays known as \emph{classical selection principles} is going
back to the papers \cite{menger}, \cite{hurewicz1},
\cite{hurewicz2}, \cite{rothberger38}. In this paper we shall
briefly discuss these classical selection principles and their
relations with other fields of mathematics, and after that we
shall be concentrated on recent innovations in the field,
preferably on results which are not included into two nice survey
papers by M. Scheepers \cite{marionnis}, \cite{marionlecce}. In
particular, in \cite{marionlecce} some information regarding
"modern", non-classical selection principles can be found. No
proofs are included in the paper.

Two classical selection principles are defined in the following
way.

Let ${\mathcal A}$ and ${\mathcal B}$ be sets whose elements are
families of subsets of an infinite set $X$. Then:

\smallskip
\noindent $\sfin({\mathcal A},{\mathcal B})$ denotes the selection
hypothesis:
\begin{quote}
For each sequence $(A_n:n\in\naturals)$ of elements of ${\mathcal
A}$ there is a sequence $(B_n:n\in\naturals)$ of finite sets such
that for each $n$, $B_n\subset A_n$, and
$\bigcup_{n\in\naturals}B_n$ is an element of ${\mathcal B}$.
\end{quote}

\smallskip \noindent  $\sone({\mathcal A},{\mathcal B})$ denotes
the selection principle:
\begin{quote}
For each sequence $(A_n:n\in\naturals)$ of elements of ${\mathcal
A}$ there is a sequence $(b_n:n\in\naturals)$ such that for each
$n$, $b_n\in A_n$, and $\{b_n:n\in\naturals\}$ is an element of
${\mathcal B}$.
\end{quote}

\medskip
In \cite{menger} Menger introduced a property for metric spaces
(called now the \emph{Menger basis property}) and in
\cite{hurewicz1} Hurewicz has proved that that property is
equivalent to the property $\sfin({\mathcal O},{\mathcal O})$,
where $\mathcal O$ denotes the family of open covers of the space,
and called now the \emph{Menger property}.

In the same paper (see also \cite{hurewicz2}) Hurewicz introduced
another property, nowadays called the \emph{Hurewicz property},
defined as follows.  A space $X$ has the Hurewicz property if for
each sequence $({\mathcal U}_n:n\in \naturals)$ of open covers
there is a sequence $({\mathcal V}_n:n\in \naturals)$ such that
for each $n$, ${\mathcal V}_n$ is a finite subset of ${\mathcal
U}_n$ and each element of the space belongs to all but finitely
many of the sets $\cup{\mathcal V}_n$. It was shown in \cite{coc7}
that the Hurewicz property can be expressed in terms of a
selection principle of the form $\sfin({\mathcal A},{\mathcal
B})$.

A selection principle of the $\sone(\mathcal A, \mathcal B)$ type
was introduced in the 1938 Rothberger's paper \cite{rothberger38},
in connection with his study of strong measure zero sets in metric
spaces that were first defined by Borel in \cite{borel}. The
\emph{Rothberger property} is the property $\sone(\mathcal O,
\mathcal O)$.

Other two properties of this sort were introduced by Gerlits and
Nagy in \cite{gerlitsnagy} under the names \emph{$\gamma$-sets}
and \emph{property $(*)$} (the later is of the $\sone(\mathcal A,
\mathcal B)$ kind as it was shown in \cite{coc7}).

\medskip
The collections $\mathcal{A}$ and $\mathcal{B}$ that we consider
here will be mainly families of open covers of some topological
space. We give now the definitions of open covers which are
important for this exposition.

\smallskip
An open cover $\mathcal U$ of a space $X$ is:
\begin{itemize}
\item an \emph{$\omega$-cover} if $X$ does not belong to $\mathcal
U$ and every finite subset of $X$ is contained in a member of
$\mathcal U$ \cite{gerlitsnagy}.
\item a \emph{$k$-cover} if $X$ does not belong to $\mathcal
U$ and every compact subset of $X$ is contained in a member of
$\mathcal U$ \cite{mccoy}.
\item a \emph{$\gamma$-cover}  if it is infinite and each $x\in X$
belongs to all but finitely many elements of $\mathcal U$
\cite{gerlitsnagy}.
\item a \emph{$\gamma_k$-cover} if each compact subset of $X$ is
contained in all but finitely many elements of $\mathcal U$ and
$X$ is not a member of the cover \cite{gamahyp}.
\item \emph{large} if each $x\in X$ belongs to infinitely many
elements of $\mathcal U$ \cite{coc1}.
\item \emph{groupable} if it can be expressed as a countable
union of finite, pairwise disjoint subfamilies $\mathcal U_n$,
$n\in\naturals$, such that each $x\in X$ belongs to $\cup\mathcal
U_n$ for all but finitely many $n$ \cite{coc7}.
\item \emph{$\omega$-groupable} if it is an $\omega$-cover and is
a countable union of finite, pairwise disjoint subfamilies
$\mathcal U_n$, $n\in\naturals$, such that each finite subset of
$X$ is contained in some element $U$ in $\mathcal U_n$ for all but
finitely many $n$ \cite{coc7}.
\item \emph{weakly groupable} if  it is a countable union of
finite, pairwise disjoint sets $\mathcal U_n$, $n\in\naturals$,
such that for each finite set $F\subset X$ we have $F\subset \cup
\mathcal U_n$ for some $n$ \cite{coc8}.
\item a \emph{$\tau$-cover} if it is large and for any two distinct
points $x$ and $y$ in $X$ either the set $\{U\in\mathcal{U}:x\in U
\mbox{ and } y\notin U\}$ is finite, or the set
$\{U\in\mathcal{U}:y\in U \mbox{ and } x\notin U\}$ is finite
\cite{tsaban-tau}.
\item a \emph{$\tau^*$-cover} if it is large and for each $x$
there is an infinite set $A_x\subset\{U\in\mathcal{U}:x\in U\}$
such that whenever $x$ and $y$ are distinct, then either
$A_x\setminus A_y$ is finite, or $A_y\setminus A_x$ is finite
\cite{st2}.
\end{itemize}

\noindent For a topological space $X$ we denote:
\begin{itemize}
\item $\Omega$ -- the family of $\omega$-covers of $X$;
\item $\mathcal K$ -- the family of $k$-covers of $X$;
\item $\Gamma$ -- the family of $\gamma$-covers of $X$;
\item $\Gamma_k$ -- the family of $\gamma_k$-covers of $X$;
\item $\Lambda$ -- the family of large covers of $X$;
\item $\mathcal O\sp{gp}$ -- the family of groupable covers of $X$;
\item $\Lambda\sp{gp}$ -- the family of groupable large covers of $X$;
\item $\Omega\sp{gp}$ -- the family of $\omega$-groupable covers of $X$;
\item $\mathcal O^{wgp}$ -- the family of weakly groupable covers of
$X$;
\item ${\sf T}$ -- the set of $\tau$-covers of $X$;
\item ${\sf T}^*$ -- the set of $\tau^*$-covers of $X$.
\end{itemize}

All covers that we consider are \emph{infinite and countable}
(spaces whose each $\omega$-cover contains a countable subset
which is an $\omega$-cover are called \emph{$\omega$-Lindelof} or
\emph{$\epsilon$-spaces} and spaces whose each $k$-cover contains
a countable subset that is a $k$-cover are called
\emph{$k$-Lindel\"of}).

So we have
\begin{center}
$\Gamma\subset {\sf T} \subset {\sf T}^* \subset \Omega \subset
\Lambda \subset \mathcal{O}$,\\
$\Gamma_k \subset \Gamma \subset \Omega\sp{gp}\subset \Omega
\subset \Lambda\sp{wgp} \subset \Lambda \subset \mathcal O$,\\
$\Gamma_k \subset \mathcal K \subset \Omega$.\\
\end{center}

\noindent In this notation, according to the definitions and
results mentioned above, we have:
\begin{itemize}
\item The Menger property: $\sfin(\mathcal O,\mathcal O)$;
\item The Rothberger property: $\sone(\mathcal O,\mathcal O)$;
\item The Hurewicz property: $\sfin(\Omega, \Lambda\sp{gp})$;
\item The $\gamma$-set property: $\sone(\Omega,\Gamma)$;
\item The Gerlits-Nagy property $(*)$: $\sone(\Omega,
\Lambda\sp{gp})$.
\end{itemize}

\noindent It is also known:
\begin{itemize}
\item $X\in \sfin(\Omega,\Omega)$ iff $ (\forall n\in\naturals) \,
X\sp n \in \sfin(\mathcal O,\mathcal O)$ \cite{coc2};
\item $X\in \sone(\Omega,\Omega)$ iff $ (\forall n\in\naturals) \,
X\sp n \in \sone(\mathcal O,\mathcal O)$ \cite{sakai1};
\item $X\in \sfin(\Omega,\Omega\sp{gp})$ iff $ (\forall n\in\naturals) \,
X\sp n \in \sfin(\Omega,\Lambda\sp{gp})$ \cite{coc7};
\item $X\in \sone(\Omega,\Gamma)$ iff $ (\forall n\in\naturals) \,
X\sp n \in \sfin(\Omega,\Gamma)$ \cite{gerlitsnagy};
\item $X\in \sone(\Omega,\Omega\sp{gp})$ iff $ (\forall n\in\naturals) \,
X\sp n \in \sone(\Omega,\Lambda\sp{gp})$ \cite{coc7}.
\end{itemize}

\noindent For a space $X$ and a point $x\in X$ the following
notation will be used:
\begin{itemize}
\item $\Omega_x$ -- the set $\{A\subset X\setminus\{x\}:
x\in\overline{A}\}$;
\item $\Sigma_x$ -- the set of all
nontrivial sequences in $X$ that converge to $x$.
\end{itemize}

\noindent A countable element $A\in\Omega_x$ is said to be
\emph{groupable} \cite{coc7}  if it can be expressed as a union of
infinitely many finite, pairwise disjoint sets $B_n$, $n\in
\naturals$, such that each neighborhood $U$ of $x$ intersects all
but finitely many sets $B_n$. We put:
\begin{itemize}
\item $\Omega_x\sp{gp}$ -- the set of groupable elements of
$\Omega_x$.
\end{itemize}

\subsection*{Games}

Already Hurewicz observed that there is a natural connection
between the Menger property and an infinitely long  game for two
players. In fact, in \cite{hurewicz1} Hurewicz implicitly proved
that the principle $\sfin(\mathcal O, \mathcal O)$ is equivalent
to a game theoretical statement (ONE does not have a winning
strategy in the game $\gfin(\mathcal O, \mathcal O)$; see the
definition below and for the proof see \cite{coc1}).

Let us define games which are naturally associated to the
selection principles $\sfin({\mathcal A},{\mathcal B})$ and
$\sone({\mathcal A},{\mathcal B})$ introduced above.

Again, ${\mathcal A}$ and ${\mathcal B}$ will be sets whose
elements are families of subsets of an infinite set $X$.

\smallskip
$\gfin({\mathcal A},{\mathcal B})$ denotes an infinitely long game
for two players, ONE and TWO, which play a round for each positive
integer. In the $n$-th round ONE chooses a set $A_n \in {\mathcal
A}$, and TWO responds by choosing a finite set $B_n \subset A_n$.
The play $(A_1,B_1,\cdots,A_n,B_n,\cdots)$ is won by TWO if
$\bigcup_ {n\in \naturals}B_n \in {\mathcal B}$; otherwise, ONE
wins.

$\gone({\mathcal A},{\mathcal B})$ denotes a similar game, but in
the $n$-th round ONE chooses a set $A_n \in {\mathcal A}$, while
TWO responds by choosing an element $b_n\in A_n$. TWO wins a play
$(A_1,b_1; \cdots; A_n,b_n;\cdots)$ if $\{b_n:n\in \naturals\} \in
\mathcal B$; otherwise, ONE wins.

It is evident that if ONE does not have a winning strategy in the
game $\gone({\mathcal A},{\mathcal B})$ (resp. $\gfin({\mathcal
A},{\mathcal B})$) then the selection hypothesis $\sone({\mathcal
A},{\mathcal B})$ (resp. $\sfin({\mathcal A},{\mathcal B})$) is
true. The converse implication need not be always true.

We shall see that a number of selection principles we mentioned
can be characterized by the corresponding game (see Table 1).

\subsection*{Ramsey theory}

Ramsey Theory is a part of combinatorial mathematics which deals
with partition symbols. In 1930, F.P. Ramsey proved the first
important partition theorems \cite{ramsey}. Nowadays there are
many "partition symbols" that have been extensively studied. We
shall consider here two partition relations (the ordinary
partition relation and the Baumgartner-Taylor partition relation)
which have nice relations with classical selection principles and
infinite game theory in topology. M. Scheepers was the first who
realized these connections (see \cite{coc1}). In \cite{coc7} very
general results of this sort were given. They show how to derive
Ramsey-theoretical results from game-theoretic statements, and how
selection hypotheses can be derived from Ramseyan partition
relations. For a detail exposition on applications of Ramsey
theory to topological properties see \cite{isgda}.

We shall also list several results which demonstrate how some
closure properties of function spaces can be also described Ramsey
theoretically and game theoretically (see Table 2).

Let us mention that no Ramseyan results are known for
non-classical selection principles that appeared in the literature
in recent years.

We are going now to define the two partition relations we shall do
with.

For a set $X$ the symbol $[X]^n$ denotes the set of $n$-element
subsets of $X$, while ${\mathcal A}$ and ${\mathcal B}$ are as in
the definitions of selection hypotheses and games. Let $n$ and $m$
be positive integers. Then:

\medskip
\noindent The \emph{ordinary partition symbol} (or \emph{ordinary
partition relation})
\[
{\mathcal A} \rightarrow ({\mathcal B})^n_m
\]
denotes the statement:
\begin{quote}
For each $A\in{\mathcal A}$ and for each function
$f:[A]^n\rightarrow\{1,\cdots,m\}$ there are a set $B\in{\mathcal
B}$ with $B\subset A$ and some $i\in\{1,\cdots,m\}$ such that for
each $Y \in [B]^n$, $f(Y) = i$.
\end{quote}

\smallskip
\noindent  The \emph{ Baumgartner-Taylor partition symbol}
\cite{BT}
\[
{\mathcal A}\rightarrow \lceil{\mathcal B} \rceil^2_m
\]
denotes the following statement:
\begin{quote}
For each $A$ in ${\mathcal A}$ and for each function $f:[A]^2
\rightarrow \{1,\cdots,m\}$ there are a set $B\in{\mathcal B}$
with $B\subset A$, an $i\in\{1,\cdots,m\}$ and a partition $B =
\bigcup_{n\in \naturals}B_n$ of $B$ into pairwise disjoint finite
sets such that for each $\{x,y\}\in [B]^2$ for which $x$ and $y$
are not from the same $B_n$, we have $f(\{x,y\}) = i$.
\end{quote}

Several selection principles of the form $\sone({\mathcal
A},{\mathcal B})$ (resp. $\sfin({\mathcal A},{\mathcal B})$) can
be characterized by the ordinary (resp. the Baumgartner-Taylor)
partition relation (see Tables 1 and 2).

\medskip
Our topological notation and terminology are standard and follow
those from \cite{engelking} with one exception: in Section 3 and
Section 4 Lindel\"of spaces are not supposed to be regular. All
spaces are assumed to be \emph{ infinite and Hausdorff}. (Notice
that some of results which will be mentioned here hold for wider
classes of spaces than it is indicated in our statements.) For a
Tychonoff space $X$ $\cpx$ (resp. $\ckx$) denotes the space of all
continuous real-valued functions on $X$ with the topology of
pointwise convergence (resp. the compact-open topology). $\zero$
denotes the constantly zero function from $\cpx$ and $\ckx$. Some
notions we are doing with will be defined when they become
necessary.

\medskip
The paper is organized in the following way. In Section 2 we give
results showing relationships between selection principles, game
theory and partition relations, as well as showing duality between
covering properties of a space $X$ and function spaces $\cpx$ and
$\ckx$ over $X$. Section 3 is devoted to duality between covering
properties of a space $X$ (expressed in terms of selection
principles) and properties of hyperspaces over $X$ that appeared
recently in the literature. In Section 4 we discuss star selection
principles -- an innovation in selection principles theory. In
particular, we discuss selection principles in uniform spaces and
topological groups. Finally, Section 5 contains some results
concerning another innovation in the field -- relative selection
principles. Several open problems are included in each section.
For a detail exposition about open problems we refer the
interested reader to \cite{spmproblem}.

\section{Selection principles, games, partition relations}

Relationships of classical selection principles with the
corresponding games and  partition relations are given in the
following table. "Game" means "ONE has no winning strategy", $n$
and $m$ are positive integers and "Source" gives papers in which
results were originally shown - the first for games and the second
for partition relations.

\medskip In Table 2 we give some results concerning relations
between covering properties of a Tychonoff space $X$ and closure
properties of the function space $\cpx$ over $X$.

\medskip
Let us recall that a space $X$ has \emph{countable tightness}
(resp. the \emph{Fr\'echet-Urysohn property} FU) if for each $x\in
X$ and each $A\in\Omega_x$ there is a countable set $B\subset A$
with $B\in\Omega_x$ (resp. a sequence $(x_n:n\in\naturals)$ in $A$
converging to $x$). $X$ is SFU (\emph{strictly FU}) if it
satisfies $\sone(\Omega_x,\Sigma_x)$ for each $x\in X$. $X$ has
\emph{countable fan tightness} \cite{arh1} (resp. \emph{countable
strong fan tightness} \cite{sakai1}) if it satisfies
$\sfin(\Omega_x, \Omega_x)$ (resp. $\sone(\Omega_x,\Omega_x)$) for
each $x\in X$. $X$ has the \emph{Reznichenko property} (E.
Reznichenko, 1996) \cite{rezfunction}, \cite{coc7} if for every
$x\in X$ each $A\in\Omega_x$ contains a countable set $B\subset A$
with $B\in\Omega_x\sp{gp}$.

\begin{center}
\begin{tabular}{|| c | l | l | l | l ||} \hline\hline
   & Selection Property          & Game     & Partition relation  & Source \\ \hline\hline
1  & $\sone(\mathcal{O},\mathcal{O})$      &
$\gone(\mathcal{O},\mathcal{O})$    & $\Omega
\rightarrow(\Lambda)^2_m$ & \cite{pawlikowski}, \cite{coc1} \\
\hline 2 & $\sone(\Omega,\Omega)$ &
$\gone(\Omega,\Omega)$ & $\Omega\rightarrow(\Omega)^n_m$ & \cite{coc3}, \cite{coc1}\\
\hline 3 & $\sone(\mathfrak B, \mathfrak B)$ & $\gone(\mathfrak B,
\mathfrak B)$ & $\mathfrak B_\Omega \rightarrow (\mathfrak B)^2_m$ & \cite{borelcovers}\\
\hline 4 & $\sone(\mathcal K,\mathcal K)$  & ? & $\mathcal
K\rightarrow(\mathcal K)^2_2$ & \cite{k-covers}\\
\hline 5 & $\sone(\Omega,\Gamma)$
& $\gone(\Omega,\Gamma)$ & $\Omega\rightarrow(\Gamma)^2_m$ & \cite{gerlitsnagy}, \cite{coc1}\\
\hline 6 & $\sone(\mathcal K,\Gamma)$ & $\gone(\mathcal
K,\Gamma)$ & $\mathcal K \rightarrow (\Gamma)^n_m$ & \cite{k-covers}, \cite{gamram}\\
\hline 7 & $\sone(\mathcal K,\Gamma_k)$ & $\gone(\mathcal
K,\Gamma_k)$ & $\mathcal K \rightarrow (\Gamma_k)^n_m$ & \cite{gamahyp}, \cite{gamram}\\
\hline 8 & $\sone(\Omega,\Lambda^{gp})$ &
$\gone(\Omega,\Lambda^{gp})$ & $\Omega
\rightarrow(\Lambda^{gp})^n_m$ & \cite{coc7}\\
\hline 9 & $\sone(\Omega,\Omega^{gp})$ &
$\gone(\Omega,\Omega^{gp})$ & $\Omega
\rightarrow(\Omega^{gp})^n_m$ & \cite{coc7}\\
\hline 10  & $\sone(\Omega,\mathcal{O}^{wgp})$ &
$\gone(\Omega,\mathcal{O}^{wgp})$  &
$\Omega\rightarrow(\mathcal{O}^{wgp})^2_m$ & \cite{coc8}\\
\hline 11  & $\sfin(\mathcal{O},\mathcal{O})$ &
$\gfin(\mathcal{O}, \mathcal{O})$      & $\Omega \rightarrow
\lceil\mathcal{O}\rceil^2_m$  &  \cite{hurewicz1}, \cite{coc1} \\
\hline 12  & $\sfin(\Omega,\Omega)$ & $\gfin (\Omega,\Omega)$ & $
\Omega \rightarrow\lceil\Omega \rceil^2_m$ & \cite{coc3},
\cite{coc1}+ \cite{coc2}
\\ \hline 13  & $\sfin(\mathcal K,\mathcal K)$
& ? &  $\mathcal K \rightarrow\lceil\mathcal K\rceil^2_2$ & \cite{k-covers} \\
\hline 14  & $\sfin(\Omega,\Lambda^{gp})$ &
$\gfin(\Omega,\Lambda^{gp})$ &
$\Omega\rightarrow\lceil\Lambda^{gp}\rceil^2_m$ & \cite{coc7}
\\ \hline 15  & $\sfin(\Omega,\Omega^{gp})$ &
$\gfin(\Omega,\Omega^{gp})$  &
$\Omega\rightarrow\lceil\Omega^{gp}\rceil^2_m$ & \cite{coc7} \\
\hline 16  & $\sfin(\Omega,\Lambda^{wgp})$  &
$\gfin(\Omega,\Lambda^{wgp})$ &
$\Omega\rightarrow\lceil\Lambda^{wgp}\rceil^2_m$ & \cite{coc8}\\
\hline 17  & $\sfin(\Omega,({\sf T}^*)\sp{gp})$  &
$\gfin(\Omega,({\sf T}^*)\sp{gp})$ &
$\Omega\rightarrow\lceil ({\sf T}^*)\sp{gp}\rceil^2_m$ & \cite{st2}\\
\hline 18  & $\sone(\Omega,({\sf T}^*)\sp{gp})$  &
$\gone(\Omega,({\sf T}^*)\sp{gp})$ &
$\Omega\rightarrow (({\sf T}^*)\sp{gp})^2_m$ & \cite{st2}\\
\hline
\end{tabular}\\
\medskip
Table 1
\end{center}

\begin{center}
\begin{tabular}{|| c | l | l | l | l ||} \hline\hline
   & $\cpx$     & $X$     & $X^n \, (\forall n\in\naturals)$ & Source \\ \hline\hline

1 & Countable tightness &  $\omega$-Lindel\"of & Lindel\"of & \cite{arhbook} \\
\hline

2 & FU  & $\sone(\Omega,\Gamma)$ & $\sone(\Omega,\Gamma)$ & \cite{gerlitsnagy}\\
\hline

3 & SFU  & $\sone(\Omega,\Gamma)$ & $\sone(\Omega,\Gamma)$ & \cite{gerlitsnagy}\\
\hline

4 & $\sfin(\Omega_{{\zero}},\Omega_{{\zero}})$ &
$\sfin(\Omega,\Omega)$ & $\sfin(\mathcal{O},\mathcal{O})$ &
\cite{arh1}\\ \hline

5 & $\gfin(\Omega_{{\zero}},\Omega_{{\zero}})$  &
$\sfin(\Omega,\Omega)$ & $\sfin(\mathcal{O},\mathcal{O})$ &
\cite{coc3} \\ \hline

6 & $\Omega_{\zero}\rightarrow\lceil \Omega_{\zero}\rceil^2_2$ &
$\sfin(\Omega,\Omega)$ & $\sfin(\mathcal{O},\mathcal{O})$ &
\cite{coc3} \\ \hline

7 & $\sone(\Omega_{{\zero}},\Omega_{{\zero}})$  &
$\sone(\Omega,\Omega)$ & $\sone(\mathcal{O},\mathcal{O})$ &
\cite{sakai1} \\ \hline

8 & $\gone(\Omega_{{\zero}},\Omega_{{\zero}})$  &
$\sone(\Omega,\Omega)$ & $\sone(\mathcal{O},\mathcal{O})$ &
\cite{coc3} \\ \hline

9 & $\Omega_{\zero} \rightarrow (\Omega_{\zero})^n_m$ &
$\sone(\Omega,\Omega)$ & $\sone(\mathcal{O},\mathcal{O})$ &
\cite{coc3} \\ \hline

10 & $\sfin(\Omega_{\zero},\Omega_{\zero}^{gp})$ &
$\sfin(\Omega,\Omega\sp{gp})$ & $\sfin(\Omega,\Lambda^{gp})$ & \cite{coc7}\\
\hline

11 & $\gfin(\Omega_{\zero},\Omega_{\zero}^{gp})$ &
$\sfin(\Omega,\Omega\sp{gp})$ & $\sfin(\Omega,\Lambda^{gp})$ & \cite{coc7}\\
\hline

12 & $\Omega_{\zero}\rightarrow\lceil
\Omega_{\zero}^{gp}\rceil^2_m$ &
$\sfin(\Omega,\Omega\sp{gp})$ & $\sfin(\Omega,\Lambda^{gp})$ & \cite{coc7}\\
\hline

13 & $\sone(\Omega_{\zero},\Omega_{\zero}^{gp})$ &
$\sone(\Omega,\Omega\sp{gp})$ & $\sone(\Omega,\Lambda^{gp})$ & \cite{coc7}\\
\hline

14 & $\gone(\Omega_{\zero},\Omega_{\zero}^{gp})$ &
$\sone(\Omega,\Omega\sp{gp})$ & $\sone(\Omega,\Lambda^{gp})$ & \cite{coc7}\\
\hline

15 & $\Omega_{\zero} \rightarrow (\Omega_{\zero}^{gp})^n_m$ &
$\sone(\Omega,\Omega\sp{gp})$ & $\sone(\Omega,\Lambda^{gp})$ & \cite{coc7}\\
\hline
\end{tabular}\\
\medskip Table 2
\end{center}

\medskip
The item 9 in Table 2 says that each finite power of a Tychonoff
space $X$ has the Hurewicz property if and only if $\cpx$ has
countable fan tightness as well as the Reznichenko property, while
the item 12 states that all finite powers of $X$ have the
Gerlits-Nagy property $(*)$ if and only if $\cpx$ has countable
strong fan tightness and Reznichenko's property. In
\cite{sakairezpyt}, conditions under which $\cpx$ has \emph{only}
the Reznichenko property have been found.

An $\omega$-cover $\mathcal U$ is \emph{$\omega$-shrinkable} if
for each $U\in\mathcal U$ there exists a closed set $C(U) \subset
U$ such that $\{C(U):U\in\mathcal U\}$ is a closed $\omega$-cover
of $X$.

\begin{theorem} [\cite{sakairezpyt}] \label{cpxrez} For a Tychonoff
space $X$ the space $\cpx$ has the Rezni\-chenko property if and
only if for each $\omega$-shrinkable $\omega$-cover is
$\omega$-groupable.
\end{theorem}

 In \cite{sakaiwfu} it was
shown that for each analytic space $X$ the space $\cpx$ has the
Reznichenko property.

Let us also mention that some other closure
properties of $\cpx$ spaces can be characterized by covering
properties of $X$ ($T$-tightness and set-tightness
\cite{sakaivariations}, selective strictly $A$-space property in
\cite{vlada}, the Pytkeev property in \cite{sakairezpyt}).

\medskip
Some results regarding the function space $\ckx$ are listed in the
following table.

\begin{center}
\begin{tabular}{|| c | l | l | l ||} \hline\hline
   & $\ckx$     & $X$   & Source \\ \hline\hline

1 & Countable tightness & $k$-Lindel\"of & \cite{mccoy}, \cite{okuyamaterada} \\
\hline

2 & SFU  & $\sone(\mathcal K,\Gamma_k)$ & \cite{chinese}\\
\hline

3 & $\sone(\Omega_{{\zero}},\Omega_{{\zero}})$  & $\sone(\mathcal
K, \mathcal K)$ &  \cite{ckx}
\\ \hline

4 & $\sfin(\Omega_{{\zero}},\Omega_{{\zero}})$ &
$\sfin(\mathcal K,\mathcal K)$ & \cite{chinese}\\
\hline

\end{tabular}
\end{center}

\begin{center}
Table 3
\end{center}

\medskip
We mention also the following result from \cite{ckx}:

\begin{theorem} \label{rezckx} If  $ONE$ has no winning strategy in the
game $\gfin(\mathcal K, \mathcal K\sp {gp})$ $($resp.
$\gone(\mathcal{K}, \mathcal{K}^{gp}))$ on $X$, then $\ckx$ has
the Reznichenko property and countable fan tightness $($resp.
countable strong fan tightness$)$.
\end{theorem}

Similarly to Theorem \ref{cpxrez} one proves

\begin{theorem} \label{ckxrez} For a Tychonoff
space $X$ the space $\ckx$ has the Reznichenko property if and
only if for each $k$-shrinkable $k$-cover is $k$-groupable.
\end{theorem}

Clearly, we say that a  $k$-cover $\mathcal U$ of a space is
\emph{$k$-shrinkable} if for each $U\in\mathcal U$ there exists a
closed set $C(U) \subset U$ such that $\{C(U):U\in\mathcal U\}$ is
a closed $k$-cover of the space.

\medskip
Quite recently it was shown how some results from Table 2 can be
applied to get pure topological characterizations of several
classical covering properties in terms of continuous images into
the space $\reals\sp\omega$ (Ko\v cinac).

\medskip
We end this section by some open problems; some of them seem to be
difficult.

The next four problems we borrow from \cite{coc8} (the first two
of them were formulated in \cite{coc2} in a different form).

\begin{problem}
Is $\sfin(\Gamma,\Lambda^{wgp}) = \sfin(\Gamma,\Omega)$?
\end{problem}

\begin{problem}
If the answer to the previous problem is "not", does
$\sfin(\Gamma,\Lambda^{gp})$ imply $\sfin(\Gamma,\Omega)$?
\end{problem}

\begin{problem}
Is $\sone(\Omega,\Lambda^{wgp})$ stronger than
$\sone(\Omega,\Lambda)$?
\end{problem}

\begin{problem}
Is $\sone(\Omega,\Omega)$ stronger than
$\sone(\Omega,\Lambda^{wgp})$?
\end{problem}

A set $X$ of reals is said to be a $\tau$-set (Tsaban) if each
$\omega$-cover of $X$ contains a countable family which is a
$\tau$-cover?

The next two problems are taken from \cite{spmproblem}.

\begin{problem}
Is the $\tau$-set property equivalent to the $\gamma$-set
property?
\end{problem}

\begin{problem}
Is $\sfin(\Omega,{\sf T})$ equivalent to the $\tau$-set property?
\end{problem}

According to \cite{k-covers} (resp. \cite{gamahyp}) a space $X$ is
a \emph{$k$-$\gamma$-set} (resp. \emph{$\gamma_k$-set};
\emph{$\gamma_k^\prime$-set}) if it satisfies the selection
hypotheses $\sone(\mathcal K,\Gamma)$ (resp. each $k$-cover
$\mathcal U$ of $X$ contains a countable set
$\{U_n:n\in\naturals\}$ which is a $\gamma_k$-cover; satisfies
$\sone(\mathcal K,\Gamma_k)$). From these two papers we take the
next two problems.

\begin{problem} Is the $k$-$\gamma$-set property equivalent to the
assertion that each $k$-cover of $X$ contains a sequence which is
a $\gamma$-cover?
\end{problem}

\begin{problem} Is the the $\gamma_k$-set property equivalent to
$\sone(\mathcal K,\Gamma_k)$?
\end{problem}

\begin{problem}
Is the converse in Theorem $\ref{rezckx}$ true?
\end{problem}


\section{Hyperspaces}

In this section we discuss duality between properties of a space
$X$ and spaces of closed subsets of $X$ with different topologies
illustrating how a selection principle for $X$ can be described by
properties of a hyperspace over $X$. As we shall see, this duality
often looks as duality between $X$ and function spaces over $X$.

By $2^X$ we denote the family of all closed subsets of a space $X$.
For a subset $A$ of $X$ we put
\begin{center}
$A^c=X \setminus A$, \, \,   $A^{+} =\{F\in 2^X:F\subset A\}$, \, \,
$A^{-} =\{F\in 2^X:F\cap A\neq\emptyset\}$.
\end{center}

The most known and popular among topologies on $2^X$ is the
Vietoris topology ${\sf V} = {\sf V^-} \vee {\sf V^+}$, where the
\emph{lower Vietoris topology} ${\sf V}^-$ is generated by all
sets $A^-$, $A\subset X$ open, and the \emph{upper Vietoris
topology} ${\sf V}^+$ is generated by sets $B^+$, $B$ open in $X$.

However, we are interested in other topologies on $2^X$.

Let $\Delta$ be a subset of $2^X$. Then the \emph{upper
$\Delta$-topology}, denoted by $\Delta^+$ \cite{poppe} is the
topology whose subbase is the collection
\[
\{(D^c)^+:D \in \Delta\} \cup \{2^X\}.
\]
Note: if $\Delta$ is closed for finite unions and contains all
singletons, then the previous collection is a base for the
$\Delta^+$-topology. We consider here two important special cases:

1. $\Delta$ is the family of all finite subsets of $X$, and

2. $\Delta$ is the collection of compact subsets of $X$.

\noindent The corresponding
$\Delta^+$-topologies will be denoted by ${\sf Z}^+$ and ${\sf
F^+}$, respectively and both have the collections of the above
kind as basic sets.  The ${\sf F}^+$-topology is known as the
\emph{upper Fell topology} (or the \emph{co-compact topology})
\cite{fell}. The \emph{Fell topology} is ${\sf F} =
{\sf \Delta^+} \vee {\sf V^-}$, where ${\sf V^-}$ is the lower
Vietoris topology.

\medskip
A number of results concerning selection principles in hyperspaces
with the $\Delta^+$-topologies obtained in the last years is listed
in the following two tables. We would like to say that some results
going to a similar direction can be found in \cite{chv}, \cite{hou}.

\begin{center}
\begin{tabular}{|| c | l | l | l ||} \hline\hline
   &  $\zplusx$ & $(\forall Y)$($Y$ open in $X$)
     & Source \\ \hline\hline
1  & countable tightness  & $\omega$-Lindel\"of & folklore \\ \hline

2  & FU  & $\sone(\Omega,\Gamma)$ & \cite{gamahyp} \\ \hline

3  & SFU & $\sone(\Omega,\Gamma)$
& \cite{gamahyp} \\ \hline

4  & countable fan tightness  & $\sfin(\Omega,\Omega)$
& \cite{selhyp} \\ \hline

5  & countable strong fan tightness  & $\sone(\Omega,\Omega)$
& \cite{selhyp} \\ \hline

6  & $(\forall S\in 2^X)$ $\sfin(\Omega_S, \Omega_S^{gp})$  &
$\sfin(\Omega,\Omega^{gp})$ & \cite{rezpythyp} \\ \hline

7  & $(\forall S\in 2^X)$ $\sone(\Omega_S, \Omega_S^{gp})$  &
$\sone(\Omega,\Omega^{gp})$ & \cite{rezpythyp} \\ \hline\hline
\end{tabular}
\end{center}

\begin{center}
Table 4
\end{center}

\begin{center}
\begin{tabular}{|| c | l | l | l ||} \hline\hline
   &  $\fplusx$ & $(\forall Y)$($Y$ open in $X$)
     & Source \\ \hline\hline

1  & countable tightness  & $k$-Lindel\"of & \cite{chv} \\ \hline

2  & FU  & $\gamma_k$-set & \cite{gamahyp} \\ \hline

3  & SFU & $\sone(\mathcal K,\Gamma_k)$
& \cite{gamahyp} \\ \hline

4 & countable fan tightness  & $\sfin(\mathcal K,\mathcal K)$
& \cite{selhyp} \\ \hline

5  & countable strong fan tightness  & $\sone(\mathcal K,\mathcal
K)$ & \cite{selhyp} \\ \hline

6  & $(\forall S\in 2^X)$ $\sfin(\Omega_S, \Omega_S^{gp})$  &
$\sfin(\mathcal K,\mathcal K^{gp})$ & \cite{rezpythyp} \\ \hline

7  & $(\forall S\in 2^X)$ $\sone(\Omega_S, \Omega_S^{gp})$  &
$\sone(\mathcal K,\mathcal K^{gp})$ & \cite{rezpythyp} \\ \hline

\end{tabular}
\end{center}

\begin{center}
Table 5
\end{center}

\medskip

One more nice property of a space has been considered in a number of
recent papers.

\smallskip
Call a space $X$ \emph{selectively Pytkeev} \cite{rezpythyp} if
for each $x\in X$ and each sequence $(A_n:n\in\naturals)$ of sets
in $\Omega_x$  there  is an infinite family
$\{B_n:n\in\naturals\}$ of countable infinite sets which is a
$\pi$-network at $x$ and such that for each $n$, $B_n\subset A_n$.
If all the sets $A_n$ are equal to a set $A$, one obtains the
notion of \emph{Pytkeev spaces} introduced in \cite{pytkeev} and
then studied in \cite{malyhintironi} (where the name Pytkeev space
was used), \cite{fedeliledonne}, \cite{sakairezpyt},
\cite{rezpythyp}. It was shown in \cite{rezpythyp} that (from a
more general result) we have the following.

\begin{theorem}\label{pytkeevgeneral} For a space $X$ the following
are equivalent:
\begin{itemize}
\item[$(1)$] $\fplusx$ has the selectively Pytkeev property;
\item[$(2)$] For each open set $Y \subset X$ and each sequence
$(\mathcal U_n:n\in\naturals)$ of $k$-covers of $Y$ there is a
sequence $(\mathcal V_n:n\in\naturals)$ of infinite, countable
sets such that for each $n$, $\mathcal V_n\subset \mathcal U_n$
and $\{\cap\, \mathcal V_n:n \in\naturals\}$ is a (not necessarily
open) $k$-cover of $Y$.
\end{itemize}
\end{theorem}

Similar assertions (from \cite{rezpythyp}) can be easily formulated
for the selectively Pytkeev property in $\zplusx$ and the Pytkeev
property in both $\zplusx$ and $\fplusx$.

\smallskip
Every (sub)sequential space has the Pytkeev property \cite[Lemma
2] {pytkeev} and every Pytkeev space has the Reznichenko property
\cite[Corollary 1.2]{malyhintironi}.

It is natural to ask.

\begin{problem}
If $\fplusx$ has the Pytkeev property, is $\fplusx$ sequential?
What about $\zplusx$?
\end{problem}

However, for a locally compact Hausdorff spaces $X$ the
countable tightness property, the Reznichenko property and the
Pytkeev property coincide in the space $(2^X,{\sf F})$ and each of
them is equivalent to the fact that $X$ is both hereditarily
separable and hereditarily Lindel\"of. There are some models of
${\sf ZFC}$ in which each of these properties is equivalent
to sequentiality of $(2^X, {\sf F})$ (for locally compact Hausdorff
spaces) \cite{rezpythyp}.

\medskip
At the end of this section we shall discuss the
Arhangel'ski\v i
$\alpha_i$ properties \cite{arhalpha} of hyperspaces according to
\cite{alphahyp}.

A space $X$ has property:
\begin{itemize}
\item[$\alpha_1$:] if for each $x\in X$ and each sequence $(\sigma_n:
n \in \naturals)$ of elements of $\Sigma_x$ there is a $\sigma \in
\Sigma_x$ such that for each $n\in\naturals$ the set $\sigma_n
\setminus \sigma$ is finite;
\item[$\alpha_2$:] if for each $x\in X$ and each sequence $(\sigma_n:
n \in \naturals)$ of elements of $\Sigma_x$ there is a $\sigma \in
\Sigma_x$ such that for each $n\in\naturals$ the set $\sigma_n
\cap \sigma$ is infinite;
\item[$\alpha_3$:] if for each $x\in X$ and each sequence $(\sigma_n:
n \in \naturals)$ of elements of $\Sigma_x$ there is a $\sigma \in
\Sigma_x$ such that for infinitely many $n\in\naturals$ the set
$\sigma_n \cap \sigma$ is infinite;
\item[$\alpha_4$:] if for each $x\in X$ and each sequence $(\sigma_n:
n\in\naturals)$ of elements of $\Sigma_x$ there is a $\sigma \in
\Sigma_x$ such that for infinitely many $n\in\naturals$ the set
$\sigma_n \cap \sigma$ is nonempty.
\end{itemize}
It is understood that
\[
\alpha_1 \Rightarrow \alpha_2 \Rightarrow \alpha_3 \Rightarrow
\alpha_4.
\]

In \cite{alphahyp} it is shown a result regarding the
$\Delta^+$-topologies one of whose corollaries is the following
theorem.

\begin{theorem} \label{alpha234} For a space $X$ the following
statements are equivalent:
\begin{itemize}
\item[$(1)$] $\fplusx$ is an $\alpha_2$-space;
\item[$(2)$] $\fplusx$ is an $\alpha_3$-space;
\item[$(3)$] $\fplusx$ is an $\alpha_4$-space;
\item[$(4)$] For each $S\in 2\sp X$, \, $\fplusx$ satisfies
$\sone(\Sigma_S,\Sigma_S)$;
\item[$(5)$] Each open set $Y \subset X$ satisfies $\sone(\Gamma_k,
\Gamma_k)$.
\end{itemize}
\end{theorem}

At the very end of this section we emphasize the existence of
results that have been appeared in
the literature \cite{daniels}, \cite{coc5}, \cite{qaPR} in
connection with selection principles in the \emph{Pixley-Roy space}
${\sf PR}(X)$ over $X$ -- the set of finite subsets of $X$ with the
topology whose base form the sets
\[
[F,U]:= \{S\in {\sf PR}(X): F\subset S\subset U\},
\]
where $F$ is a finite and $U$ is an open set in $X$ with $F\subset
U$.

\section{Star and uniform selection principles}

We repeat that in this section we assume that all topological spaces
are Hausdorff and $\omega$-Lindel\"of.

In \cite{starmenger}, Ko\v cinac introduced star selection
principles in the following way.

\medskip
Let ${\mathcal A}$ and ${\mathcal B}$ be collections of open
covers of a space $X$ and let $\mathcal K$ be a family of subsets
of $X$. Then:

\noindent {\bf 1.} The symbol $\ssone(\mathcal A,\mathcal B)$ denotes the
selection hypothesis:
\begin{quote}
For each sequence $(\mathcal U_n:n\in\naturals)$ of elements of
$\mathcal A$ there exists a sequence $(U_n:n\in\naturals)$ such
that for each $n$, $U_n\in\mathcal U_n$ and $\{$St$(U_n,\mathcal
U_n):n\in\naturals\}$ is an element of $\mathcal B$;
\end{quote}

\noindent {\bf 2.} The symbol $\ssfin(\mathcal A,\mathcal B)$ denotes the
selection hypothesis:
\begin{quote}
For each sequence $(\mathcal U_n:n\in\naturals)$ of elements of
$\mathcal A$ there is a sequence $(\mathcal V_n:n \in \naturals)$
such that for each $n\in \naturals$, $\mathcal V_n$ is a finite
subset of $\mathcal U_n$, and $\bigcup_{n \in \naturals}
\{$St$(V,\mathcal U_n):V \in\mathcal V_n\} \in \mathcal B$;
\end{quote}

\noindent {\bf 3.} By $\sufin(\mathcal A,\mathcal B)$ we denote the
selection hypothesis:
\begin{quote}
For each sequence $(\mathcal U_n:n \in\naturals)$ of members of
$\mathcal A$ there exists a sequence $(\mathcal V_n:n \in
\naturals)$ such that for each $n$, $\mathcal V_n$ is a finite
subset of $\mathcal U_n$ and $\{$St$(\cup \mathcal V_n,\mathcal
U_n):n \in \naturals\}\in \mathcal B$.
\end{quote}

\noindent {\bf 4.} ${\sf SS}_{\mathcal K}^{*}(\mathcal A,\mathcal B)$
denotes the selection hypothesis:
\begin{quote}
For each sequence $(\mathcal U_n:n\in \naturals)$ of elements of
$\mathcal A$ there exists a sequence $(K_n:n \in \naturals)$ of
elements of $\mathcal K$ such that $\{$St$(K_n, \mathcal U_n):n
\in \naturals\} \in\mathcal B$.
\end{quote}
\noindent When $\mathcal K$ is the collection of all one-point (resp.,
finite) subspaces of $X$ we write $\sssone(\mathcal
A,\mathcal B)$ (resp., $\sssfin(\mathcal A,\mathcal B)$)
instead of ${\sf SS}_{\mathcal K}^{*}(\mathcal A,\mathcal B)$.

\smallskip
Here, for a subset $A$ of a space $X$ and a collection
$\mathcal S$ of subsets of $X$, \, St$(A,\mathcal S)$ denotes the
star of $A$ with respect to $\mathcal S$, that is the set $\cup\{S
\in \mathcal S:A \cap S \neq \emptyset \}$; for $A =\{x\}$, $x \in
X$, we write St$(x, \mathcal S)$ instead of $St(\{x\},\mathcal
S)$.

The following terminology we borrow from \cite{starmenger}.
\smallskip
\noindent
A space $X$ is said to have:

1.  the \emph{star-Rothberger property} SR,

2.  the \emph{star-Menger property} SM,

3.  the \emph{strongly star-Rothberger property} SSR,

4. the {\it strongly star-Menger property} SSM,

\noindent if it satisfies the selection hypothesis:

1.  $\ssone(\mathcal O,\mathcal O)$,

2.  $\ssfin(\mathcal O,\mathcal O)$ (or, equivalently, $\sufin(\mathcal O,
\mathcal O)$),

3. $\sssone(\mathcal O,\mathcal O)$,

4. $\sssfin(\mathcal O,\mathcal O)$,

\noindent respectively.

In \cite{starhur}, two star versions of the Hurewicz property were
introduced as follows:
\begin{itemize}
\item[SH:]{A space $X$ satisfies the \emph{star-Hurewicz property} if
for each sequence $({\mathcal U}_n:n\in\naturals)$ of open covers
of $X$ there is a sequence $({\mathcal V}_n:n\in\naturals)$ such
that for each $n\in\naturals$ \, ${\mathcal V}_n$ is a finite
subset of ${\mathcal U} _n$ and each $x\in X$ belongs to
St$(\cup{\mathcal V}_n,\mathcal U_n)$ for all but finitely many
$n$.}
\item[SSH:]{A space $X$ satisfies the \emph{strongly star-Hurewicz
property} if for each sequence $({\mathcal U}_n:n\in\naturals)$ of
open covers of $X$ there is a sequence $(A_n:n\in\naturals)$ of
finite subsets of $X$ such that each $x\in X$ belongs to
St$(A_n,\mathcal U_n)$ for all but finitely many $n$ (i.e. if $X$
satisfies $\sssfin(\mathcal O,\Gamma)$).}
\end{itemize}

Of course Menger spaces are SSM, amd every SSM space is SM.
Similarly for the Hurewicz and Rothberger properties.

There is a strongly star-Menger space which is not Menger, but
every metacompact strongly star-Menger space is a Menger space
\cite{starmenger}. For paracompact (Hausdorff) spaces the three
properties, SM, SSM and M, are equivalent \cite{starmenger}. The
same situation is with the classes SSH, SH and H \cite{starhur}.

The product of two star-Menger (resp. SH) spaces need not be in
the same class. But if one factor is compact, then the product is
in the same class \cite{starmenger}, \cite{starhur}. A Lindel\"of
space is not a preserving factor for classes SSM and SSH.

In \cite{starmenger} we posed the following still open problem.

\begin{problem} Characterize spaces $X$ which are
SM (SSM, SR, SSR) in all finite powers.
\end{problem}

A partial solution of this problem was given in \cite{starhur}.

\begin{theorem}\label{smpowers} If each finite power of a space $X$
is $SM$, then $X$ satisfies $\sufin(\mathcal O,\Omega)$.
\end{theorem}

\begin{theorem}\label{ssmpowers} If all finite powers of a space
$X$ are strongly star-Menger, then $X$ satisfies $\sssfin(\mathcal
O,\Omega)$.
\end{theorem}

In the same paper we read the following two assertions.

\begin{theorem}\label{ufin-wgp} For a space $X$ the following are equivalent:
\begin{itemize}
\item[$(1)$] $X$ satisfies $\sufin(\mathcal O,\Omega)$;
\item[$(2)$] $X$ satisfies $\sufin(\mathcal O,\mathcal O^{wgp})$.
\end{itemize}
\end{theorem}

\begin{theorem}\label{ssfin-wgr} For a space $X$ the following are equivalent:
\begin{enumerate}
\item[$(1)$] {$X$ satisfies $\sssfin(\mathcal O,\Omega);$}
\item[$(2)$]{$X$ satisfies $\sssfin(\mathcal O,\mathcal O^{wgp})$.}
\end{enumerate}
\end{theorem}

So the previous problem can be now translated to

\begin{problem}
Does $X \in \sufin(\mathcal O,\mathcal O^{wgp})$ imply that all
finite powers of $X$ are star-Menger? Is it true that
$\ssfin(\mathcal O,\Omega) = \ssfin(\mathcal O,\mathcal O^{wgp})$?
Does $X \in \sssfin(\mathcal
O,\mathcal O^{wgp})$ imply that each finite power of $X$ is SSM?
\end{problem}

The following result regarding star-Hurewicz spaces

\begin{theorem}\label{sshgroup} For a space $X$ the following are
equivalent:
\begin{enumerate}
\item[$(1)$]{$X$ has the strongly star-Hurewicz property;}
\item[$(2)$]{$X$ satisfies the selection principle $\sssfin(\mathcal O,
\mathcal O^{gp})$.}
\end{enumerate}
\end{theorem}

\noindent suggests the following

\begin{problem} Is it true that  $\ssfin(\mathcal O,\Gamma)=
\ssfin(\mathcal O,\mathcal O^{gp})$?
\end{problem}

Let us formulate once again a question from \cite{starmenger}.

\begin{problem} Characterize hereditarily SM (SSM, SR, SSR, SH, SSH)
spaces.
\end{problem}

Let $X$ be a space. Two players, ONE and TWO, play a round per
each natural number $n$. In the $n$--th round ONE chooses an open
cover $\mathcal U_n$ of $X$ and TWO responds by choosing a finite
set $A_n\subset X$. A play $\mathcal U_1,A_1;\cdots;\mathcal
U_n,A_n;\cdots$ is won by TWO if $\{$St$(A_n,\mathcal U_n):n \in
\naturals\}$ is a $\gamma$-cover of $X$; otherwise, ONE wins.

\smallskip
Evidently, if ONE has no winning strategy in the strongly
star-Hurewicz game, then $X$ is an SSH space.

\begin{conjecture}
The strongly star-Hurewicz property of a space $X$ need not imply
ONE does not have a winning strategy in the strongly star-Hurewicz
game played on $X$.
\end{conjecture}

Similar situation might be expected in cases of star versions of
the Menger and Rothberger properties and the corresponding games
(which can be naturally associated to a selection principle).

But the situation can be quite different in case of
zero-dimensional metrizable topological groups (see the next
section).

\bigskip
In \cite{selunif} it was demonstrated that selection principles in
uniform spaces are a good application of star selection principles
to concrete special classes of spaces. In particular case of
topological groups ones obtain nice classes of groups.

Recall that a uniformity on a set $X$ can be defined in terms of
uniform covers, and then the uniform space is viewed as the pair
$(X,\mathbb C)$, or in terms of entourages of the diagonal, and
then the uniform space is viewed as the pair $(X,\mathbb U)$
\cite{engelking}. The first approach is convenient because it
allows us to define uniform selection principles in a natural way
similar to the definitions of topological selection principles.
After that it is easy to pass to  $(X,\mathbb U)$.

Let us explain this on the example of the uniform Menger property.
A uniform space $(X,\mathbb C)$ is \emph{uniformly Menger} or
\emph{Menger-bounded} if for each sequence $(\alpha_n:n\in
\naturals)$ of uniform covers there is a sequence
$(\beta_n:n\in\naturals)$ of finite
sets such that for each $n\in\naturals$, $\beta_n\subset \alpha_n$ and
$\bigcup_{n\in\naturals}\beta_n$ is a (not necessarily uniform)
cover of $X$.

\begin{theorem} \label{menger}For a uniform space $(X,\mathbb C)$
the following are equivalent:
\begin{itemize}
\item[$(a)$] {$X$ has the uniform Menger property;}
\item[$(b)$] {for each sequence $(\alpha_n:n\in \naturals)\subset
\mathbb C$ there is a sequence $(A_n:n\in\naturals)$ of finite
subsets of $X$ such that
$X=\bigcup_{n\in\naturals}$St$(A_n,\alpha_n)$;}
\item[$(c)$] {for each sequence $(\alpha_n:n\in \naturals)\subset
\mathbb C$ there is a sequence $(\beta_n:n\in\naturals)$ such that
for each $n$ \, $\beta_n$ is a finite subset of $\alpha_n$ and $X=
\bigcup_{n\in\naturals}$St$(\cup\beta_n,\alpha_n)$.}
\end{itemize}
\end{theorem}

Therefore, we conclude that here we have, in notation we adopted,
that $\sfin(\mathbb C,\mathcal O) = \sssfin(\mathbb C,\mathcal O)=
\ssfin(\mathbb C,\mathcal O)$. In other words, one can say that a
uniform space $(X,\mathbb U)$ is uniformly Menger if and only if
for each sequence $(U_n:n\in\naturals)$ of entourages of the
diagonal of $X$ there is a sequence $(A_n:n\in\naturals)$ of
finite subsets of $X$ such that
$X=\bigcup_{n\in\naturals}U_n[A_n]$.

\medskip
It is understood, if a uniform space $X$ has the Menger property
with respect to topology generated by the uniformity, then $X$ is
uniformly Menger. However, any non-Lindel\"of Tychonoff space
serves as an example of a space which is uniformly Menger that has
no the Menger property. (Similar remarks hold for the uniform
Rothberger and uniform Hurewicz properties defined below.) But a
regular topological space $X$ has the Menger property if and only
if its fine uniformity has the uniform Menger property. Uniform
spaces having the uniform Menger property have some properties
which are similar to the corresponding properties of totally
bounded uniform spaces.

\medskip
In case of topological groups we have: A topological group
$(G,\cdot)$ is \emph{Menger-bounded} if for each sequence
$(U_n:n\in\naturals)$ of neighborhoods of the neutral element
$e\in G$ there is a sequence $(A_n:n\in\naturals)$ of finite
subsets of $G$ such that $X=\bigcup_{n\in\naturals}A_n\cdot U_n$.
This class of groups was already studied in the literature under
the name $o$-bounded groups \cite{hernandez}, \cite{tkachenko}.

More information on Menger-bounded topological groups the reader
can find in \cite{hernandez}, \cite{tkachenko}, \cite{coc11},
\cite{taras1}, \cite{taras2}, \cite{tsabangroups}.

\medskip
Similarly, a uniform space $(X,\mathbb C)$ is \emph{Rothberger-bounded}
if it satisfies one of the three equivalent selection hypotheses:
$\sone(\mathbb C,\mathcal O)$, $\sssone(\mathbb C,\mathcal O)$,
$\ssone(\mathbb C,\mathcal O)$.

A topological group $(G,\cdot)$
is \emph{Rothberger-bounded} if for each sequence $(U_n:n\in\naturals)$ of
neighborhoods of the neutral element $e\in G$ there is a sequence
$(x_n:n\in\naturals)$ of elements of $G$ such that
$X=\bigcup_{n\in\naturals}x_n\cdot U_n$.

\medskip
Finally, a uniform space $(X, \mathbb C)$
is  \emph{uniformly Hurewicz} if for each sequence $(\alpha_n:n\in
\naturals)$ of uniform covers of $X$ there is a sequence
$(F_n:n\in\naturals)$ of finite subsets of $X$ such that each
$x\in X$ belongs to all but finitely many sets  St$(F_n,\alpha_n)$.

It is easy to define \emph{Hurewicz-bounded} topological groups.

The difference between uniform and topological selection
principles is big enough \cite{selunif}. Here we point out some of
differences on the example of the Hurewicz properties. (Note that
uniformly Hurewicz spaces have many similarities with totally
bounded uniform spaces.)

Every subspace of a uniformly Hurewicz uniform space
is uniformly Hurewicz.
A uniform space $X$ is uniformly Hurewicz if and
only if its completion $\tilde{X}$ is uniformly Hurewicz.
The product of two uniformly Hurewicz uniform spaces
is also uniformly Hurewicz.

Hurewicz-bounded topological groups are preserving factors for the
class of Menger-bounded groups \cite{coc11}.

\section{Relative selection principles}

A systematic study of relative topological properties was started
by A.V. Arhangel'ski\v i in 1989 and then continued in a series of
his papers and papers of many other authors (see for example
\cite{arhrel1}, \cite{arhrel2}).

Let $X$ be a topological space and $Y$ a subspace of $X$. To each
topological property $\mathcal P$ (of $X$) associate a property
"relative $\mathcal P$" which shows how $Y$ is located in $X$;
thus we speak also that \emph{$Y$ is relatively $\mathcal P$ in
$X$}. For $Y=X$ the relative version of a property $\mathcal P$
must be just $\mathcal P$. In that sense classical topological
properties are called \emph{absolute} properties.

A systematic investigation of relative selection principles was
initiated by Ko\v cinac (see, \cite{fareast},
\cite{kocguidobabin}, \cite{guidokoc}). Later on it was shown that
relative covering properties described by selection principles,
like absolute ones, have nice relations with game theory and
Ramsey theory, as well as with with measure-like and basis-like
properties in metric spaces and topological groups. We shall see
that relative selection principles can be quite different from
absolute ones. For example, in \cite{coc8} it was shown that a
very strong relative covering property is not related to a weak
absolute covering property. More precisely, it was proved that the
Continuum Hypothesis implies the existence of a relative
$\gamma$-subset $X$ of the real line such that $X$ does not have
the (absolute) Menger property $\sfin(\mathcal{O},\mathcal{O})$.
It will be also demonstrated that relative selection principles
strongly depend on the nature of the basic space.

Notice that much still needs to be investigated regarding the
relative selection principles in connection with "non-classical"
selection principles.

 Let $X$ be a space and $Y$ a subset of $X$. We use the symbol
$\mathcal O_X$ to denote the family of open covers of $X$ and the
symbol $\mathcal O_Y$ for the set of covers of $Y$ by sets open in
$X$. Similar notation will be used for other families of covers.

\noindent In this notation we have:

\begin{itemize}
\item \emph{$Y$ is relatively Menger in $X$} \cite{fareast}
\item \emph{$Y$ is relatively Rothberger in $X$} \cite{fareast}
\item \emph{$Y$ is relatively Hurewicz in $X$} \cite{guidokoc},
\cite{coc8}
\item \emph{$Y$ is a relative $\gamma$-set in $X$} \cite{kocguidobabin}
\end{itemize}
\noindent if the following selection principle is satisfied
\begin{itemize}
\item{$\sfin(\mathcal O_X,\mathcal O_Y)$}
\item{$\sone(\mathcal O_X,\mathcal O_Y)$}
\item{$\sfin(\Omega_X,\mathcal O_Y\sp{gp})$}
\item{$\sone(\Omega_X,\Gamma_Y)$.}
\end{itemize}

\noindent When $Y=X$ we obtain considered absolute versions of
selection principles.

In Section 2 we saw that there is a nice duality between covering
properties of a Tychonoff space $X$ and closure properties of
function spaces $\cpx$ and $\ckx$. In what follows we show that
similar duality exists between relative selection principles and
closure properties of mappings. For a Tychonoff space $X$ and its
subspace $Y$ the restriction mapping $\pi:\cpx \rightarrow \cpy$
is defined by $\pi(f) = f\restriction Y$, $f\in \cpx$.

\subsection*{Relative Menger property}

If $f:X\rightarrow Y$ is a continuous mapping, then we say that
$f$ has \emph{countable fan tightness} if for each $x\in X$ and
each sequence $(A_n:n\in \naturals)$ of elements of $\Omega_x$
there is a sequence $(B_n:n \in\naturals)$ of finite sets such
that for each $n$, $B_n\subset A_n$ and $\bigcup_{n\in \naturals}
f(B_n) \in \Omega_{f(x)}$.

\smallskip
The following theorem from \cite{fareast} gives a relation between
relative Menger-like properties and fan tightness of mappings.

\begin{theorem}\label{relMengerpi} For a Tychonoff space $X$
and a subspace $Y$ of $X$ the following are equivalent:
\begin{itemize}
\item[$(1)$] For all $n\in \naturals$, $Y^n$ is Menger in $X^n$;
\item[$(2)$]  $\sfin(\Omega_X,\Omega_Y)$ holds;
\item[$(3)$] The mapping $\pi$ has countable fan tightness.
\end{itemize}
\end{theorem}

In \cite{liljanateza} the relative Menger property was further
considered and the following theorem proved (compare with item 11
in Table 1):

\begin{theorem}\label{relmengergame} Let $X$ be a Lindel\"of space.
Then for each subspace $Y$ of $X$ the following are equivalent:
\begin{enumerate}
\item[$(1)$]{$\sfin(\mathcal{O}_X,\mathcal{O}_{Y})$;}
\item[$(2)$]{ONE has no winning strategy in $\gfin(\mathcal{O}_X,\mathcal{O}_{Y})$;}
\item[$(3)$]{For each natural number $m$, $\Omega_{X} \rightarrow \lceil
\mathcal{O}_{Y}\rceil^2_m$.}
\end{enumerate}
\end{theorem}

The following result from \cite{coc10} is of the same sort and is
a relative version of a result from \cite{coc8}.

\begin{theorem}\label{relfinwgp} Let $X$ be a space with
the the Menger property $\sfin(\mathcal{O}_X,\mathcal{O}_X)$ and
$Y$ a subspace of $X$. The following are equivalent:
\begin{enumerate}
\item[$(1)$]{$\sfin(\Omega_X,\mathcal{O}^{wgp}_{Y})$;}
\item[$(2)$]{ONE has no winning strategy in the game
$\gfin(\Omega_X,\mathcal{O}^{wgp}_{Y})$;}
\item[$(3)$]{For each $m\in\naturals$, $\Omega_X\rightarrow\lceil
\mathcal{O}^{wgp}_{Y}\rceil^2_m$.}
\end{enumerate}
\end{theorem}

The relative Menger property in metric spaces has basis-like and
measure-like characterizations as it was shown in \cite{coc9} and
\cite{coc10}. Relative Menger-like properties in topological
groups also have very nice characterizations \cite{coc11}.

To formulate results in this connection we need some terminology.

In \cite{menger} Menger introduced a property for metric spaces
$(X,d)$ that we call  the \emph{Menger basis property}: For each
base $\mathcal{B}$ in $X$ there is a sequence $(B_n:n\in
\naturals)$ in $\mathcal{B}$ such that $\lim_{n\rightarrow\infty}
diam(B_n) = 0$ and the set $\{B_n:n\in \naturals\}$ is an open
cover of $X$. As we mentioned in Introduction, in \cite{hurewicz1}
W. Hurewicz proved that a metrizable space $X$ has the Menger
basis property with respect to all metrics on $X$ generating the
topology of $X$ if and only if it has the Menger property
$\sfin(\mathcal{O},\mathcal{O})$.

Say that a subspace $Y$ of a metric space $(X,d)$ has the
\emph{Menger basis property in $X$} if for each base $\mathcal{B}$
in $X$ there is a sequence $(B_n:n\in\naturals)$ in $\mathcal{B}$
such that $\lim_{n\rightarrow\infty} diam(B_n) = 0$ and the set
$\{B_n:n\in \naturals\}$ is an open cover of $Y$.

The following definition is motivated by the definition of strong
measure zero sets introduced by Borel in \cite{borel} (see the
subsection on relative Rothberger property).

A metric space $(X,d)$ has \emph{Menger measure zero} if for each
sequence $(\epsilon_n:n\in\naturals)$ of positive real numbers
there is a sequence $({\mathcal V}_n:n\in\naturals)$ such that:
\begin{itemize}
\item[$(i)$] {for each $n$, ${\mathcal V}_n$ is a finite family of subsets of $X$;}
\item[$(ii)$] {for each $n$ and each $V\in {\mathcal V}_n$, diam$_d(V) < \epsilon_n$;}
\item[$(iii)$]{$\bigcup_{n\in\naturals}{\mathcal V}_n$ is an open cover of $X$.}
\end{itemize}

Combining some results from \cite{coc9} and \cite{coc10} we have
the following theorem.

\begin{theorem}\label{relmengerbasis-measure} Let $(X,d)$ be a separable
zero-dimensional metric space and let $Y$ be a subspace of $X$.
The following statements are equivalent:
\begin{itemize}
\item[$(1)$] {$Y$ is relatively Menger in $X$;}
\item[$(2)$] {$Y$ has the Menger basis property in $X$;}
\item[$(3)$] {$Y$ has Menger measure zero with respect to each
metric on $X$ which gives $X$ the same topology as $d$.}
\end{itemize}
\end{theorem}

Let $(G,\cdot)$ be a topological group and $H$ its subgroup.
Denote by ${\sf M}(G,H)$ the following game for two players, ONE
and TWO, which play a round for each $n\in \naturals$. In the
$n$-th round ONE chooses a neighborhood $U_n$ of the neutral
element of $G$ and then TWO chooses a finite set $F_n\subset G$.
Two wins a play $U_1,F_1; U_2,F_2; ...$ if and only if $\{F_n\cdot
U_n:n \in \naturals\}$ covers $H$. (It is a relative version of a
game first mentioned in \cite{tkachenko}.)

In \cite{coc11}, the following result regarding Menger-like
properties for topological groups has been obtained.

\begin{theorem}\label{mengergroups} Let $G$ be a zero-dimensional
metrizable group and let $H$ be a subgroup of $G$. The following
assertions are equivalent:
\begin{itemize}
\item[$(1)$] {$H$ is Menger-bounded;}
\item[$(2)$]{$H$ is Menger-bounded in $G$;}
\item[$(3)$]{ONE has no winning strategy in the game ${\sf
M}(H,H)$;}
\item[$(4)$] $H$ has the relative Menger property in $G$;
\item[$(5)$] $H$ has Menger measure zero with respect to all
metrizations of $G$.
\end{itemize}
\end{theorem}

\subsection*{Relative Hurewicz property}

Recall that a subspace $Y$ of a space $X$ is relatively Hurewicz
in $X$ if the selection principle $\sfin(\Omega_X,\mathcal
O_Y\sp{gp}$) holds.

Following \cite{guidokoc} and \cite{coc8} we say that a continuous
mapping $f:X\rightarrow Y$ has the \emph{ selectively Reznichenko
property} if for each sequence $(A_n:n\in\naturals)$ from
$\Omega_x$ there is a sequence $(B_n:n\in\naturals)$ such that for
each $n$, \, $B_n$ is a finite subset of $A_n$ and
$\bigcup_{n\in\naturals}B_n \in\Omega_{f(x)}\sp{gp}$.

The theorem below is a combination of results from \cite{guidokoc}
and \cite{coc7} and  gives a characterization of the relative
Hurewicz property in all finite powers \cite{guidokoc}.

\begin{theorem}\label{relativeHurmapping} For a Tychonoff
space $X$ and its subspace $Y$ the following are equivalent:
\begin{itemize}
\item[$(1)$] $\pi$ has the selectively Reznichenko property;
\item[$(2)$] For each $n\in\naturals$, $Y^n$ has the Hurewicz property
in $X^n$;
\item[$(3)$]{ONE has no winning strategy in $\gfin(\Omega_X,\Omega^{gp}_{Y})$;}
\item[$(4)$]{For each $m\in\naturals$, $\Omega_{X} \rightarrow \lceil
\Omega^{gp}_{Y}\rceil^2_m$.}
\end{itemize}
\end{theorem}

The relative Hurewicz property has also a game-theoretic and
Ramsey-theoretic description \cite{coc8}.

\begin{theorem}\label{relhurewicz} Let $X$ be a Lindel\"of space.
Then for each subspace $Y$ of $X$ the following are equivalent:
\begin{enumerate}
\item[$(1)$]{$\sfin(\Omega_X,\mathcal{O}^{gp}_{Y})$;}
\item[$(2)$]{ONE has no winning strategy in $\gfin(\Omega_X,\mathcal{O}^{gp}_{Y})$;}
\item[$(3)$]{For each $m\in\naturals$, $\Omega_{X} \rightarrow
\lceil \mathcal{O}^{gp}_{Y}\rceil^2_m$.}
\end{enumerate}
\end{theorem}

Following \cite{coc8}, we are going now to show that the relative
Hurewicz property for metric spaces can be characterized by
basis-like and measure-like properties.

Let $(X,d)$ be a metric space and $Y$ a subspace of $X$. Then:

$Y$ has the \emph{Hurewicz basis property} in $X$ if for any basis
$\mathcal{B}$ of $X$ there is a sequence $(U_n:n\in\naturals)$ in
$\mathcal{B}$ such that $\{U_n:n\in\naturals\}$ is a groupable
cover of $Y$ and $\lim_{n\rightarrow\infty} diam_d(U_n) = 0$.

$Y$  has \emph{Hurewicz measure zero} (in $X$) if for each
sequence $(\epsilon_n:n\in\naturals)$ of positive real numbers
there is a sequence $({\mathcal V}_n:n\in\naturals)$ such that:
\begin{itemize}
\item[$(i)$] {for each $n$, ${\mathcal V}_n$ is a finite family of subsets of $X$;}
\item[$(ii)$] {for each $n$ and each $V\in {\mathcal V}_n$, diam$_d(V) < \epsilon_n$;}
\item[$(iii)$]{$\bigcup_{n\in\naturals}{\mathcal V}_n$ is a groupable cover of $X$.}
\end{itemize}

\begin{theorem}[\cite{coc8}]\label{relHurbasis-measure} Let $(X,d)$ be
a metric space and let $Y$ be a subspace of $X$. The following
statements are equivalent:
\begin{itemize}
\item[$(1)$] {$Y$ is relatively Hurewicz in $X$;}
\item[$(2)$] {$Y$ has the Hurewicz basis property in $X$.}
\end{itemize}
If $(X,d)$ is zero-dimensional and separable, then conditions
$(1)$ and $(2)$ are equivalent to
\begin{itemize}
\item[$(3)$] {$Y$ has Hurewicz measure zero with respect to each
metric on $X$ which gives $X$ the same topology as $d$ does.}
\end{itemize}
\end{theorem}

For special topological groups we have interesting
characterizations of relative versions of Hurewisz-like properties
\cite{coc11}. The following result shows again how relative
properties depend on the structure of the basic space.

\begin{theorem}\label{hurewiczgroups}
For a subgroup $(G,+)$ of $(^{\omega}\integers,+)$ the following
are equivalent:
\begin{itemize}
\item[$(1)$] $G$ is Hurewicz-bounded;
\item[$(2)$] $G$ has Hurewicz measure zero in the Baire metric on
$^{\omega}\integers$;
\item[$(3)$] $G$ has the relative Hurewicz property in
$^{\omega}\integers$.
\end{itemize}
\end{theorem}

Similar results for the selection principle $\sfin(\Omega_X,
\mathcal O_Y\sp{wgp})$ can be found in \cite{liljanateza},
\cite{coc9} and \cite{coc10}.

\subsection*{Relative Rothberger property}

A continuous mapping $f:X\rightarrow Y$ is said to have
\emph{countable strong fan tightness} if for each $x\in X$ and
each sequence $(A_n:n\in \naturals)$ from $\Omega_x$ there exist
$x_n\in\ A_n$, $n\in \naturals$, such that $\{f(x_n):n\in
\naturals\} \in \Omega_{f(x)}$.

\smallskip
Here is a theorem from \cite{fareast} which gives a
characterization of the relative Rothberger property.

\begin{theorem}\label{relRothbergerpi} If $Y$ is a subset of a
Tychonoff space $X$ then the following are equivalent:
\begin{itemize}
\item[$(1)$] For each $n\in \naturals$, $Y^n$ has the Rothberger
property in $X^n$;
\item[$(b)$] The selection principle $\sone(\Omega_X,\Omega_Y)$ holds;
\item[$(c)$] The mapping $\pi:\cpx \rightarrow \cpy$ has countable strong fan tightness.
\end{itemize}
\end{theorem}

In \cite{borel} Borel defined a notion for metric spaces $(X,d)$
nowadays called \emph{strong measure zero}.  $Y\subset X$ is
\emph{strong measure zero} if for each sequence
$(\epsilon_n:n\in\naturals)$ of positive real numbers there is a
sequence $(U_n:n\in\naturals)$ of subsets of $X$ such that for
each $n$, $U_n$ has diameter $< \epsilon_n$ and
$\{U_n:n\in\naturals\}$ covers $Y$.

In \cite{MF} it was shown that a metric space $X$ has the
(absolute) Rothberger property if and only if it has strong
measure zero with respect to each metric on $X$ which generates
the topology of $X$.

A result from \cite{marionsmz} states that if $Y$ is a subset of a
$\sigma$-compact metrizable space $X$, then $Y$ has the relative
Rothberger property in $X$ if and only if $Y$ has strong measure
zero with respect to each metric on $X$ which generates the
topology of $X$. It is also shown that in this case the previous
two conditions are equivalent to each of the next two conditions:
\begin{itemize}
\item{ONE has no winning strategy in the game $\gone(\Omega_X,\mathcal{O}_{Y})$;}
\item{For each $m\in\naturals$, $\Omega_X\rightarrow(\mathcal{O}_{Y})^2_m$.}
\end{itemize}

To state a result similar to (a part of) Theorem
\ref{relmengerbasis-measure} we need the following notion.

A subset $Y$ of a metric space $(X,d)$ has the \emph{Rothberger
basis property} in $X$ if for each base $\mathcal{B}$ in $X$ and
for each sequence $(\epsilon_n:n\in\naturals)$ of positive real
numbers there is a sequence $(B_n:n\in \naturals)$ of elements of
$\mathcal{B}$ such that  $diam_d(B_n) < \epsilon_n$, and
$\{B_n:n\in\naturals\}$ covers $Y$.

\begin{theorem}\label{rothbasis} Let $X$ be a metrizable space, $Y$
a subspace of $X$. Then  the following are equivalent:
\begin{itemize}
\item[$(1)$]$Y$ has the relative Rothberger property in $X$;
\item[$(2)$]$Y$ has the Rothberger basis property in $X$ with respect
to all metrics generating the topology of $X$.
\end{itemize}
\end{theorem}

For Rothberger-bounded subgroups of the set of real numbers we
have the following description \cite{coc11}:

\begin{theorem}\label{rothbergergroups}
For a subgroup $(G,+)$ of $(\reals,+)$ the following are
equivalent:
\begin{itemize}
\item[$(1)$] $G$ is Rothberger-bounded;
\item[$(2)$] $G$ has strong measure zero in $\reals$;
\item[$(3)$] $G$ has the relative Rothberger property in $\reals$.
\end{itemize}
\end{theorem}

\subsection*{Relative Gerlits-Nagy property $(*)$}

It is the property $\sone(\Omega_X,\mathcal O\sp{gp})$, where $Y$
is a subspace of a space $X$.

We state here only one statement regarding this property in metric
spaces.

A subspace $Y$ of a metric space $(X,d)$ has the
\emph{Gerlits-Nagy basis property} in $X$ if for each base
$\mathcal{B}$ for the topology of $X$ and for each sequence
$(\epsilon_n:n\in\naturals)$ of positive real numbers there is a
sequence $(B_n:n\in\naturals)$ such that for each $n$,
$B_n\in\mathcal{B}$ and $diam(B_n) < \epsilon_n$, and
$\{B_n:n\in\naturals\}$ is a groupable cover of $Y$.

The following result is from \cite{coc9} and \cite{coc10}.

\begin{theorem}\label{relgerlitsnagy} Let $X$ be an infinite
$\sigma$-compact metrizable space and let $Y$ be a subspace of
$X$. The following statements are equivalent:
\begin{itemize}
\item[$(1)$]{$\sone(\Omega_X,\mathcal{O}^{gp}_{Y})$;}
\item[$(1)$]{ONE has no winning strategy in the game $\gone(\Omega_X,\mathcal{O}^{gp}_{Y})$;}
\item[$(1)$]{For each positive integer $m$, $\Omega_X\rightarrow(\mathcal{O}^{gp}_{Y})^2_m$;}
\item[$(1)$]{$Y$ has the Gerlits-Nagy basis property in $X$ with respect to
all metrics generating the topology of $X$.}
\end{itemize}
\end{theorem}

Let us point out that similar assertion (for $\sigma$-compact
metrizable spaces) is true for the principle
$\sone(\Omega_X,\mathcal{O}^{wgp}_{Y})$.

\subsection*{Relative $\gamma$-sets}

A subspace $Y$ of a space $X$ is a relative $\gamma$-set in $X$ if
the selection hypothesis $\sone(\Omega_X,\Gamma_Y)$ holds
\cite{kocguidobabin}.

Clearly, every $\gamma$-set is also a relative $\gamma$-set, but
the converse is not true. The relative $\gamma$-set property is
hereditary, while the Gerlits-Nagy $\gamma$-set property is not
hereditary \cite{kocguidobabin}. Relative $\gamma$-sets of real
numbers have strong measure zero. By the well known facts on
strong measure zero sets (Borel's conjecture that no uncountable
set of real numbers has strong measure zero is undecidable in
ZFC), the question if there is an uncountable relative
$\gamma$-set of real numbers is undecidable in ZFC.

The relative $\gamma$-set property, as other relative covering
properties defined in terms of selection principles, depend on the
basic space $X$. Recently, A.W. Miller \cite{millerrelgama}
considered relative $\gamma$-sets in $2\sp\omega$ and $\omega\sp
\omega$. He defined two cardinals ${\mathfrak p}(2\sp\omega)$
(resp. ${\mathfrak p}(\omega\sp\omega)$) to be the smallest
cardinality of the set $X$ in $2\sp\omega$ (resp. in ${\omega
\sp\omega}$) which is not a relative $\gamma$-set in the
corresponding space, observed that ${\mathfrak p}\le {\mathfrak
p}(\omega \sp \omega) \le {\mathfrak p}(2\sp\omega)$ and proved
that it is relatively consistent with ${\sf ZFC}$ that ${\mathfrak
p} = {\mathfrak p}(\omega \sp \omega) < {\mathfrak p}(2\sp\omega)$
and ${\mathfrak p} < {\mathfrak p}(\omega\sp\omega) = {\mathfrak
p}(2\sp\omega)$.

Here ${\mathfrak p}$ is the pseudointersection number and is the
cardinality of the smallest non-gamma set (according to a result
from \cite{gerlitsnagy}; see also \cite{coc2}).

To characterize relative $\gamma$-sets in $\reals$ we need the
following notion \cite{kocguidobabin}.

A continuous mapping $f:X\rightarrow Y$ is said to be
\emph{strongly Fr\'echet} if for each $x\in X$ and each sequence
$(A_n:n\in \naturals)$ in $\Omega_x$ there is a sequence
$(B_n:n\in \naturals)$ such that for each $n$, $B_n$ is a finite
subset of $A_n$ and the sequence $(f(B_n):n\in\naturals)$
converges to $f(x)$.

\begin{theorem}[\cite{kocguidobabin}]\label{relativegamma} For a
Tychonoff space $X$ and its subspace $Y$ the following are
equivalent:
\begin{itemize}
\item[$(1)$] $Y$ is a $\gamma$-set in $X$;
\item[$(2)$] For each $n\in \naturals$, $Y^n$ is a
$\gamma$-set in $X^n$;
\item[$(3)$] The mapping $\pi:\cpx\rightarrow \cpy$ is strongly Fr\'echet.
\end{itemize}
\end{theorem}

\subsection*{Relative SSH spaces}

We close this section by a relative version of a star selection
principle from Section 4 and considered in \cite{starhur}. Once
more we conclude that relative selection principles are very
different from absolute ones.

Let $Y$ be a subspace of a space $X$. We say that \emph{$Y$ is
strongly star-Hurewicz in $X$} if for each sequence $(\mathcal
U_n:n\in\naturals)$ of open covers of $X$ there is a sequence
$(A_n:n\in\naturals)$ of finite subsets of $X$ such that each
point $y\in Y$ belongs to all but finitely many sets
St$(A_n,\mathcal U_n)$.

There is a strongly star-Menger space $X$ and a subspace $Y$ of
$X$ such that $Y$ is relatively strongly star-Hurewicz in $X$ but
not (absolutely) strongly star-Hurewicz. The space $X$ is the
Mr\'owka-Isbel space $\Psi(\mathcal A)$ \cite{engelking}, $Y$ is
the subspace $\mathcal A$, where $\mathcal A$ is  an almost
disjoint family of infinite subsets of positive integers having
cardinality $< \mathfrak b$ (see \cite{starhur}).

\smallskip
Notice that the following two relative versions of the SSH
property could be investigated.

(1) For each sequence $(\mathcal U_n:n\in\naturals)$ of covers of
$Y$ by sets open in $X$ there is a sequence $(A_n:n\in\naturals)$
of finite subsets of $X$ such that each point $y\in Y$ belongs to
all but finitely many sets St$(A_n,\mathcal U_n)$.

(2) For each sequence $(\mathcal U_n:n\in\naturals)$ of open
covers of $X$ there is a sequence $(A_n:n\in\naturals)$ of finite
subsets of $Y$ such that each point $y\in Y$ belongs to all but
finitely many sets St$(A_n,\mathcal U_n)$.

\bigskip \noindent {\bf Acknowledgement.} This survey is written
on the basis of a lecture I delivered on the \emph{Third Seminar
on Geometry and Topology} held in Tabriz, July 15-17, 2004. My
warmest thanks to the organizers of the Seminar for their kind
invitation and the hospitality during the Seminar.

\end{document}